\documentclass{amsart}
\usepackage{palatino,hhline, bbm}
\usepackage{amsthm, amsmath}
\usepackage{amssymb} 
\usepackage{mathrsfs}
\usepackage{url}
\usepackage[hmargin=2.7cm, vmargin=3cm]{geometry}
\usepackage[breaklinks,bookmarksopen,bookmarksnumbered]{hyperref}
\usepackage[all]{xy}
\usepackage{tikz-cd}

\theoremstyle{plain}
\newtheorem{thm}{Theorem}[section]
\newtheorem*{theo}{Theorem}
\newtheorem{lem}{Lemma}[section]
\newtheorem{prop}{Proposition}[section]
\newtheorem{cor}{Corollary}[section]

\theoremstyle{remark}
\newtheorem{rem}{Remark}[section]

\newcommand\pr{\noindent\textit{Proof} : }
\def\lr#1{\langle {#1} \rangle}

\newcommand\Aut{\operatorname{Aut}}

\newcommand\codim{\operatorname{codim}}

\newcommand\im{\operatorname{Im}}
\newcommand\Ker{\operatorname{Ker}}

\newcommand\divi{\operatorname{div}}
\newcommand\pp{\scriptscriptstyle\bullet}
\newcommand{\mo}{\smallsetminus}

\newcommand\Z{\mathbb{Z}}
\newcommand\A{\mathbb{A}}

\newcommand\C{\mathbb{C}}
\renewcommand\P{\mathbb{P}}

\newcommand\G{\mathbb{G}}
\renewcommand\O{\mathscr{O}}
\newcommand\sy{\mathsf{S}}

\newcommand\iso{\vbox{\hbox to .8cm{\hfill{$\scriptstyle\sim$}\hfill}
		\nointerlineskip\hbox to .8cm{{\hfill$\longrightarrow $\hfill}} }}
\newcommand\bir{\vbox{\hbox to .8cm{\hfill{$\scriptstyle\sim$}\hfill}
		\nointerlineskip\hbox to .8cm{{\hfill$\dasharrow $\hfill}} }}
\newcommand\abs[1]{\lvert {#1}\rvert}

\newcommand\oo{^{\mathrm{o}}}
\newcommand\ooo{^{\mathrm{oo}}}

\begin{document}
	\title{Symmetric tensors on the intersection of two quadrics and Lagrangian fibration}
	\author{A. Beauville}
	\address{Universit\'e C\^ote d'Azur\\
		CNRS -- Laboratoire J.-A. Dieudonn\'e\\
		Parc Valrose\\
		F-06108 Nice cedex 2, France}
	\email{arnaud.beauville@univ-cotedazur.fr}
	\author{A. Etesse}
	\address{Institut Math\'ematique de Toulouse\\
		Universit\'e Paul Sabatier, 118 Rte de Narbonne, 31400 Toulouse}
	\email{antoine.etesse@ens-lyon.fr}
	\author{A. H\"oring}
	\address{Universit\'e C\^ote d'Azur\\
		CNRS -- Laboratoire J.-A. Dieudonn\'e\\
		Parc Valrose\\
		F-06108 Nice cedex 2, France}
	\email{Andreas.Hoering@univ-cotedazur.fr}
	\author{J. Liu}
	\address{Institute of Mathematics\\
 Academy of Mathematics and Systems Science\\
 Chinese Academy of Sciences\\
 Beijing\\
 100190\\
 China}
	\email{jliu@amss.ac.cn}

	\author{C. Voisin}
	\address{Sorbonne Universit\'{e} and Universit\'{e}  Paris Cit\'e, CNRS, IMJ-PRG, F-75005 Paris, France}
	\email{{claire.voisin@imj-prg.fr}}
	\thanks{J. Liu is supported by the National Key Research and Development Program of China (No.\ 2021YFA1002300), the NSFC grants (No.\ 12001521 and No.\ 12288201) and the CAS Project for Young Scientists in Basic Research (No.\ YSBR-033).}
	\begin{abstract}
		Let  $X$  be a $n$-dimensional (smooth) intersection of two quadrics, and let  $T^*X$  be its cotangent bundle. We show that the algebra of symmetric tensors on $ X$  is a polynomial algebra in  $n$  variables. The corresponding map  $\Phi: T^*X \rightarrow  \C^n$  is a Lagrangian fibration, which admits an explicit geometric description; its general fiber is a Zariski open subset of an abelian variety,   quotient of a hyperelliptic Jacobian by a $2$-torsion subgroup. In dimension $3$,    $\Phi$ is the Hitchin fibration of the moduli space of rank $2$ bundles with fixed  determinant on a   curve of genus $ 2$. 
		
	\end{abstract}
	\maketitle 
	\section{Introduction}
	Let $X\subset \P^{n+2}_{\C}$ be a smooth $n$-dimensional complete intersection of two quadrics, with $n\geq 2$, and let $T^*X$ be its cotangent bundle. The $\C$-algebra $H^0(T^*X,\O_{T^*X})$ is canonically isomorphic to the algebra of symmetric tensors $H^0(X,\sy^{\pp}T_X )$.
	Recall that $T^*X$ carries a canonical symplectic structure. Our main result is the following theorem:
	
	\begin{theo}
		\label{theo: main}
		\emph{a)} The vector space $W:=H^0(X,\mathsf{S}^{2}T_X )$ has dimension $n$, and the natural map   $\sy^{\pp} W\rightarrow H^0(X,\sy^{\pp}T_X )$ is an isomorphism. 
		
		\emph{b)} The corresponding map $\Phi: T^*X\rightarrow W^*\cong\C^n$ is a Lagrangian fibration.
		
		\emph{c)} When $X$ is general, the general fiber of $\Phi$ is of the form $A\mo Z$, where $A$ is an abelian variety and $\codim Z\geq 2$. 
	\end{theo}
	We will  give a precise geometric description of the map $\Phi$ and of the abelian variety $A$ in \S\, \ref{fibersphi} and \ref{fibersPhi}.
	
	\subsection{Comments} 	
1)	For $n=2$, a) follows from Theorem 5.1 in \cite{DO-L}, while b) and c) are proved in \cite{K-L}. The proof is based on the isomorphism  $T_X\cong \Omega ^1_X(1)$.
	The Theorem also follows from the fact that $X$ is a moduli space for parabolic rank 2 bundles on $\P^1$ \cite{C}, so that $\Phi :T^*X\rightarrow \C^2$  is identified to the \emph{Hitchin fibration} (see \cite{B-H-K}).
	
	For $n=3$, $X$ is isomorphic to the moduli space of vector bundles of rank 2 and fixed determinant of odd degree \cite{N}; again the Theorem follows from the properties of the Hitchin fibration  (see \S \ref{n=3}).  
	It would be interesting  to have a modular interpretation of   $\Phi $ for $n \geq 4$. Note that the Hitchin map for $G$-bundles is homogeneous quadratic only when  $G$ is $\operatorname{SL}(2) $ or a product of copies of $\operatorname{SL}(2) $, so this limits the possibilities of using it.
	
	\smallskip	
	2) The map $\Phi$ is an example of an \emph{algebraically completely integrable system} --- see for instance \cite{V}, and Remark \ref{acis}. Such a
	situation is rather exceptional: most varieties do not admit nonzero symmetric tensors  (for instance,   hypersurfaces of degree $\geq 3$ \cite{H-L-S}); when they do,   even for varieties as simple as quadrics, the algebra of symmetric tensors is fairly complicated. We do not have a conceptual explanation for the particularly simple behaviour in our case.
	
	\medskip	
	3) For $n= 2$ or $3$, the generality assumption on $X$ in c) is unnecessary. It seems likely that this is the case for all $n$, but our method does not allow us to conclude.
	
	\subsection{Strategy} We will first treat the case $n=3$, which is independent of the rest of the paper (\S\,\ref{n=3}). For the general case we will develop two different approaches. In the first one we exhibit a natural $n$-dimensional subspace $W\subset H^0(X,\sy^2T_X)$, from which we deduce a map $T^*X\rightarrow W^*\cong \C^n$ (\S\,\ref{Phi}). We then show that $\Phi$ has the required properties, which implies a), b) and c) for general $X$ (\ref{res}). In the second approach (\S\,\ref{tensors}) we prove directly a) for all smooth $X$, by realizing $X$ as a double covering of a quadric.
	
	\subsection{Notations}
	Throughout the paper $X$ will be a  smooth complete intersection of two quadrics in $\P^{n+2}$, with $n\geq 2$.  We denote by $T^*X$ its cotangent bundle and by $\P T^*X$ its projectivization in the geometric sense (not in the Grothendieck sense). If $V$ is a vector space, we denote by $\P(V)$ the associated projective space $V\mo\{0\}/\C^* $ parametrising one-dimensional subspaces of $V$.

	\bigskip	 
	\section{The case $n=3$}\label{n=3}
	In this section we show how  our general results can be obtained in the case $n=3$ by interpretating $X$ as a moduli space.

	As in \ref{reid} below, we associate to $X$ a genus 2 curve $C$, such that the variety of lines in $X$ is isomorphic to $JC$.  
	Let us fix a line bundle $N$ on $C$ of degree $1$; then $X$ is  isomorphic to the moduli space $\mathscr{M}$ of rank 2 stable vector bundles on $C$ with determinant $N$ \cite{N}.
	The cotangent bundle $T^*\mathscr{M}$ is naturally identified with the moduli space of \emph{Higgs bundles}, that is pairs $(E,u)$ with $E\in\mathscr{M}$ and $u:E \rightarrow E\otimes K_C$ a homomorphism with $\operatorname{Tr}u=0 $. 
	The \emph{Hitchin map} $\Phi : T^*\mathscr{M}\rightarrow H^0(K_C^2)$ associates to a pair $(E,u)$ the section $\det u$ of $K_C^2$. It is a Lagrangian fibration \cite{H}.

	Let $\omega \in H^0(K_C^2)$.  \emph{We  assume in what follows that  $\omega $ vanishes at $4$ distinct points}. 
	Let $C_{\omega }$ be the curve in the cotangent bundle $T^*C$ defined by $z^2=\omega $. The projection  $\pi :C_{\omega }\rightarrow C$ is a double covering, branched along $\operatorname{div}(\omega ) $, and $C_{\omega }$ is a smooth curve of genus 5. Let $P$ be the Prym variety associated to $\pi $, that is, the kernel of the norm map $\operatorname{Nm}:JC_{\omega }\rightarrow JC $; it is a 3-dimensional abelian variety.
	
	\begin{prop}
		The fiber $\Phi^{-1}(\omega )$ is isomorphic to the complement of a   curve in $P$.
	\end{prop}
	\pr Recall that the map $L\mapsto \pi _*L$ establishes a bijective correspondence between line bundles on $C_{\omega }$ and rank 2 vector bundles $E$ on $C$ endowed with a homomorphism $u:E\rightarrow E\otimes K_C$ such that $u^2=\omega \cdot \operatorname{Id}_E $, or equivalently, $\operatorname{Tr}u=0 $ and $\det u=\omega $ (see for instance \cite{BNR}). 
	To get $(E,u)$ in $\Phi ^{-1}(\omega )$ we have to impose moreover $\det E= N$ and $E$ stable.
	Since $\det \pi _*L=\operatorname{Nm}(L)\otimes K_C^{-1} $, the first condition means that $L$ belongs to the translate $P_N:=\operatorname{Nm}^{-1}(K_C\otimes N) $ of $P$. 
	
	Then the vector bundle $\pi _*L$ is unstable if and only if it contains an invertible subsheaf $M$ of degree 1; this is equivalent to saying that there is a nonzero map $\pi ^*M\rightarrow L$, that is, $L=\pi ^*M(p)$ for some point $p\in C_{\omega }$. 
	The condition $L\in P_N$ means $M^{2}(\pi p)=K_C\otimes N$, so $M$ is determined by $p$ up to the 2-torsion of $JC$. Thus the locus of line bundles $L\in P_N$ such that $\pi _*L$ is unstable is a curve.\qed
	
	\medskip

	Let $\rho :C\rightarrow \P^1$ be the canonical double covering, and $B\subset\P^1$ its branch locus.  Since the homomorphism $\mathsf{S}^2H^0(K_C)\rightarrow H^0(K_C^{2})$ is surjective, the divisor of $\omega $ is of the form $\rho ^*(p+q)$, for some $p,q\in\P^1$; by assumption we have $p\neq q$ and $p,q\notin B$. 
	
	\begin{prop}
		Let $\Gamma $ be the double covering of $\P^1$ branched along $B\cup\{p,q\} $. There is an exact sequence \[\ 0\rightarrow \Z/2\rightarrow J\Gamma \rightarrow P\rightarrow 0\,.\]
	\end{prop}
	\pr Let $\chi :\P^1\rightarrow \P^1$ be  the double covering   branched along $\{p,q\} $.  Since $\operatorname{div}(\omega ) =\rho ^*(p+q)$, there is a cartesian diagram of  double coverings 
	\[\xymatrix{ C_{\omega }\ar[r]^{\xi }\ar[d]_{\pi }& \P^1\ar[d]^{\chi }\\ C\ar[r]_{\rho }& \P^1}\]
	which gives rise to two commuting involutions $\sigma ,\tau $ of $C_{\omega }$, exchanging the two sheets of $\pi $ and $\xi $ respectively. The field of rational functions on $C_{\omega }$ is   
	\[\C(x,y,z)\ |\ y^2=f(x), z^2=g(x) \]
	where $f$ and $g$ are polynomials   with $\divi f=B$ and $\divi g=\{p,q\} $. Then $\sigma $ and $\tau $ change the sign of $y$ and $z$ respectively. 
	
	The involution $\sigma \tau $ is fixed point free, so the quotient $\Gamma :=C_{\omega} /\lr{\sigma \tau} $ has genus 3; its field of functions is  $\C(x,w)$ with $w=yz$ and $w^2=f(x)g(x)$. 
	We have again a cartesian square
	\[\xymatrix{ C_{\omega }\ar[r]^{\varphi }\ar[d]_{\pi }& \,\,\Gamma \phantom{.}\ar[d]^{\psi}\\ C\ar[r]_{\rho }& \P^1.}\]
	Let $\alpha \in J\Gamma $. We have $\operatorname{Nm}_{\pi }\varphi ^*\alpha =\rho ^*\operatorname{Nm}_{\psi}\alpha =0  $, hence $\varphi ^*$ maps $J\Gamma $ into $P\subset JC_{\omega }$. Since $\varphi $ is \'etale, we have $\Ker \varphi^*=\Z/2$; since $\dim J\Gamma =\dim P=3$, $\varphi ^*$ is surjective.\qed

	\section{Definition of $\Phi$}\label{Phi}
	Let $Y$ be  a smooth degree $d$ hypersurface in $\P^N$, defined by an equation $f=0$.
	Recall that one associates to  $f$ a section $h_f $ of $\sy^2\Omega^1_Y(d) $, the \emph{hessian} or \emph{second fundamental form} of $f$ \cite{G-H}: at a point  $y$ of $Y$, the intersection of $Y$ with the tangent hyperplane  $H$ to $Y$ at $y$ is a  hypersurface in $H$ singular at $y$, and $h_{f}(y)$ is the degree 2 term in
	the Taylor expansion of $f_{|H}$ at $y$.
	
	Now let $X\subset \P^{n+r}$ be a smooth complete intersection of $r$ hypersurfaces of degree $d$; let $$V\subset H^0(\P^{n+r},\O_{\P}(d))$$ be the $r$-dimensional subspace of degree $d$ polynomials   vanishing on $X$. By restricting $h_{f}$, for $f\in V$,  to $X$, we get a linear map
	\[V\otimes \O_X\longrightarrow \sy^2\Omega^1_X(d)\]
	which gives at each point $x\in X$ a linear space of quadratic forms on the tangent space $T_{x}(X)$. Note that, when $d=2$, the corresponding quadrics in $\P (T_x(X))$ can be viewed geometrically as follows:
	the projective space $\P( T_x(X))$ can be identified with the space of lines in $\P^{n+r}$ passing through $x$ and tangent to $X$; then for each $q\in V$,  the quadric defined by $h_q(x)$ parameterizes the lines passing through $x$ and contained in the quadric $\{q=0\} $.

	Now we want to consider the ``inverse" of the quadratic form $h_f(x)$ on $T_x(X)$, that is, the form on $T^*_x(X)$ given in coordinates by the cofactor matrix.  Intrinsically,   each $f\in V$ gives a twisted symmetric morphism
	\[h_f:T_X\longrightarrow \Omega^1_X(d)\]
	which induces a twisted symmetric morphism on $(n-1)$-th exterior powers, namely
	\[\wedge^{n-1}h_f:\bigwedge\nolimits^{n-1}T_X\longrightarrow \bigwedge\nolimits^{n-1}\Omega^1_X((n-1)d)\,.\]
	We now observe that $K_X=\mathscr{O}_X(-n-1-r+dr)$, hence
	\[\bigwedge\nolimits^{n-1}T_X\cong \Omega^1_X(n+1-r(d-1))\, ,\quad\bigwedge\nolimits^{n-1}\Omega^1_X\cong T_X(-n-1+r(d-1))\, ,\] so that $\wedge^{n-1}h_f$ is in fact a symmetric morphism from
	$\Omega^1_X(n+1-r(d-1))$ to $T_X((n-1)d-n-1+r(d-1))$,
	hence provides a section $$\wedge^{n-1}h_f\in H^0(X,\sy^2T_X(d(n+2r-1)-2(n+r+1))).$$
 Being locally given by the cofactor matrix, $\wedge^{n-1}h_f$ is homogeneous of degree $n-1$ in $f$, hence we have constructed a morphism
	\[  \alpha: \sy^{n-1}V\longrightarrow H^0(X,\sy^2T_X(d(n+2r-1)-2(n+r+1))) \quad \mbox{such that }\ \alpha (f^{n-1})=\wedge^{n-1}h_f\,.\]

		From now on, we restrict to the case  $d=2,\,r=2$, so $X$ is the complete intersection of two quadrics in $\P^{n+2}$. The previous construction  gives a morphism
\[ \alpha: \sy^{n-1}V\longrightarrow H^0(X,\sy^2T_X)\,.\]

	Using the canonical isomorphism $H^0(T^*X,\O_{T^*X})=  H^0(X,\sy^{\pp}T_X)$, we deduce from $\alpha $ a morphism
	\[\Phi: T^*X\longrightarrow \sy^{n-1}V^*\cong \C^n\, .\]
	We have $\Phi(\lambda v)=\lambda ^2\Phi(v)$ for $v\in T^*X$, $\lambda \in\C$, so $\Phi$
	induces a rational map
	\[\varphi :\P T^*X \dasharrow \P^{n-1}\]
	whose indeterminacy locus $Z$ is the image of $\Phi^{-1}(0)$.
	\begin{prop}
		$1)\ \alpha $ is injective.
		
		$2)\ \Phi$ is surjective.
		
		$3)$ The image of $Z$ by the structure map $p:\P T^*X\rightarrow X$ is a proper subvariety of $X$.
	\end{prop}
	\pr Let $x$ be a general point of $X$. We claim that the base locus in $\P (T_x(X))$ of the pencil of quadratic forms $\{h_q(x)\}_{q\in V}$ is smooth. Indeed, this locus can be viewed as the variety $F_x$ of lines in $X$ passing through $x$. Let $F$ be the Fano variety of lines contained in $X$, and let \[G\subset F\times X=\{(\ell, y)\,|\,y\in\ell\}\,.\] 
	Then $F$ and therefore $G$ are smooth \cite[Theorem 2.6]{R}, hence $F_x$, which is the fiber above $x$ of the projection $G\rightarrow X$, is smooth since $x$ is general. It follows that, in an appropriate system of coordinates $(k_1,\ldots,k_n )$ of $T_x(X)$,  the forms $\{h_q(x)\}$ can be written
	\[t\sum k_i^2 +\sum \alpha _ik_i^2\quad \mbox{with}\ \alpha _i \mbox{ distinct in }\C, \ t\in \C\,.\]
	Then $\wedge^{n-1}h_q(x)$ is given by the diagonal matrix with entries $\beta _i:=\displaystyle\prod_{j\neq i}(t+\alpha _j)\ (i=1,\ldots ,n)$. These polynomials in $t$ are linearly independent, hence they generate the space of quadratic
	forms on $T^*_xX$
	which are diagonal in the basis $(k_i)$. This linear system has dimension $n$, so $\alpha $ is injective;  it has no base point, so    $\varphi $ induces a finite, surjective morphism  $\P (T^*_xX)\rightarrow \P^{n-1}$. Thus $\Phi$ is surjective, and $Z\cap \P (T^*_xX)=\varnothing$, which gives 2) and 3). \qed

	\medskip	
	
	We want to give a geometric construction of the rational map $\varphi : \P T^*X\dasharrow \P^{n-1}$. A point of $\P T^*X$ is a pair $(x,H)$, where $x\in X$ and $H$ is a hyperplane in $T_x(X)$. Restricting the pencil $\{h_q(x)\}_{q\in V} $ to $H$ gives a pencil   of quadrics on $H$, which for $(x,H)$ general contains $n-1$ singular quadrics $q_1,\ldots ,q_{n-1}$. The subset $\{q_1,\ldots ,q_{n-1}\} $ of   $\P (V)$ corresponds to a point $\varphi _{x,H}$ of $\P(\sy^{n-1}V^*)$ -- namely the hyperplane in $\sy^{n-1}V$ spanned by $q_1^{n-1},\ldots , q_{n-1}^{n-1}$. 
	
	\begin{prop}\label{geom}
		$\varphi (x,H)=\varphi _{x,H}$.
	\end{prop}
	\pr We can assume that $x$ is general. We have seen that the restriction of $\varphi $ to $\P (T^*_xX)$ is the morphism given by the linear system of quadratic forms $W\cong \sy^{n-1}V$ spanned by the forms $\wedge^{n-1}h_q(x)$, for $q\in V$; in other words, $\varphi $ maps the point $H$ of $\P (T^*_xX)$ to the hyperplane of forms in $W$ vanishing at $H$. 
	
	On the other hand, $\varphi _{x,H}$ is the hyperplane of $\sy^{n-1}V$ spanned by the $q^{n-1}$ for those $q\in V$ such that $h_q(x)_{|H}$ is singular; this condition is equivalent to say that the form $\wedge^{n-1}h_q(x)$ on $T^*_xX$ vanishes at $H$. Therefore $\varphi _{x,H}$ is spanned by quadratic forms vanishing at $H$, hence coincides with $\varphi (x,H)$.\qed
	
	\begin{cor}\label{codimZ}
		$\codim Z\geq 2$.
	\end{cor}
	\pr Suppose $Z$ contains a component $Z_0$ of codimension 1; since $p(Z)\neq X$, we have $Z_0=p^{-1}(p(Z_0))$. We claim that this is impossible, in fact $Z$ cannot contain a fiber $p^{-1}(x)$.   Indeed this would mean that for $q\in V$, the form $h_q(x)$ is singular along 
	all hyperplanes $H\subset T_xX$, that is, $h_q(x)  $ has rank $\leq n-2$.  But the rank of $h_q(x)$ is the rank of the restriction of $q$ to the projective tangent subspace to $X$ at $x$. Restricting a quadratic form to a hyperplane lowers its rank by up to two. Since a general $q$ in $V$ has rank $n+3$, its restriction to a codimension 2 subspace has rank $\geq n-1$.\qed

	\bigskip	
	\section{Fibers of $\varphi $}\label{fibersphi}
	In an appropriate system of coordinates $(x_0,\ldots ,x_{n+2})$, our variety $X$ is defined by the equations $q_1=q_2=0$, with
	\[ q_1=\sum x_i^2\ ,\quad  q_2=\sum \mu _ix_i^2\qquad\mbox{with }\,\mu _i\in \C\  \mbox{ distinct.}\]
	Let $\Pi=\P(V)\ (\cong \P^1)$ be the pencil of quadrics containing $X$. We choose a coordinate $t$ on $\Pi$ so that the quadrics of $\Pi$ are given by $tq_1-q_2=0$. Then the singular quadrics of $\Pi$ correspond to the points $\mu _0,\ldots ,\mu_{n+2} $. 
	
	The goal of this section is
	to describe  the general fiber of the rational map $\varphi :\P T^*X\dasharrow \sy^{n-1}\Pi\ (\cong \P^{n-1})$.  For $\lambda =(\lambda _1,\ldots ,\lambda _{n-1})\in \sy^{n-1}\Pi$, 
	let $C_{\mu,\lambda}$ denote the hyperelliptic curve $y^2=\prod (t-\mu _i)\prod (t-\lambda _j)$, of genus $n$. We will prove:
	
	\begin{prop}\label{birat}
		For $\lambda $ general in $\sy^{n-1}\Pi$, the fiber $\varphi ^{-1}(\lambda)$ is birational to the quotient of the Jacobian $JC_{\mu ,\lambda}$ by  the group
		$\Gamma :=\{\pm 1_{JC}\}\times \Gamma ^+ $, where $\Gamma ^+\cong  (\Z/2Z)^{n-2}$ is a group of translations by $2$-torsion elements.
	\end{prop}
	
	\subsection{Odd-dimensional intersection of 2 quadrics}\label{reid}
	We briefly recall here the results of Reid's thesis (\cite{R}, see also \cite{D-R}). Let $Y\subset \P^{2g+1}$ be a smooth intersection of 2 quadrics, and let $\Xi\ (\cong \P^1)$ be the pencil of quadrics containing $Y$. Let $\Sigma\subset \Xi $  be the subset of $2g+2$ points  corresponding to singular quadrics, and let $C$ be the double covering of $\Xi$ branched along $\Sigma $ -- this is a hyperelliptic curve of genus $g$. The intermediate Jacobian $JY$ of $Y$ is isomorphic to $JC$ (as principally polarized abelian varieties). The variety $F$ of $(g-1)$-planes contained in $Y$ is also isomorphic to $JC$, but this isomorphism is not canonical.
	
	In an appropriate system of coordinates, the equations of $Y$ are of the form
	\[\sum x_i^2=\sum \alpha _ix_i^2=0\qquad\mbox{with } \alpha _i\in \C \mbox{ distinct;}\]
	then $\Sigma =\{\alpha _1,\ldots ,\alpha _{2g+2}\} $. The group $\Gamma :=(\Z/2\Z)^{2g+1}$ acts on $Y$ (hence also on $F$)
	by changing the signs of the coordinates. Let $\Gamma ^+\subset \Gamma $ be the subgroup of elements which change an even number of coordinates. 
	For an appropriate choice of the isomorphism $F\iso JC$, 
	the image of $\Gamma^+ $ in $\Aut(JC)$ is the group $T_2$ of translations by $2$-torsion elements of $JC$, and the image of $\Gamma $ is $T_2\times \{\pm 1_{JC}\}  $ \cite[Lemma 4.5]{D-R}.

	\subsection{An auxiliary construction}\label{aux}
	We consider  the  projective space $\mathbb{P}^{2n+1}$ equipped with the  system of  homogeneous coordinates
	\[x_0,\ldots, x_{n+2};\,y_1,\ldots,\,y_{n-1}\] and the   affine space $\mathbb{A}^{n-1}$ equipped with the affine coordinates  $\lambda_1,\ldots,\lambda_{n-1}$. Let
	\[\mathscr{X}\subset  \mathbb{P}^{2n+1}\times \mathbb{A}^{n-1}\]
	be the complete intersection of the two quadrics with equations 
	\[Q_1=Q_2=0\quad\mbox{with}\quad Q_1=\sum_{i=0}^{n+2}x_i^2+\sum_{j=1}^{n-1}y_j^2\quad ,\quad 
	Q_2=\sum_{i=0}^{n+2}\mu_ix_i^2+\sum_{j=1}^{n-1}\lambda_jy_j^2\, .\]
	The second  projection $\mathscr{X}\rightarrow  \mathbb{A}^{n-1}$ gives a  family  of complete  intersections  of two  quadrics $\mathscr{X}_\lambda$ of  dimension $2n-1$ parameterized by $ \mathbb{A}^{n-1}$.  Note that $X$ is the intersection of $\mathscr{X}$ with the subspace $\P^{n+2}\subset \P^{2n+1}$ defined by $y_1=\ldots =y_{n-1}=0$.
	
	Let $p: \mathscr{F}\rightarrow \mathbb{A}^{n-1}$ be the family of  $(n-1)$-planes contained in the $\mathscr{X}_{\lambda }$, that is
	\[\mathscr{F}=\{(P,\lambda )\,|\,\lambda \in \mathbb{A}^{n-1}, \ P\ (n-1)\mbox{-plane}\subset \mathscr{X}_{\lambda }\}\,. \] 
	For $\lambda $ general, the fiber $\mathscr{F}_{\lambda }$ is isomorphic to the Jacobian of the hyperelliptic curve $C_{\mu ,\lambda}$ (\ref{reid}).

	Let $(P,\lambda )$ be a general point of $\mathscr{F}$. Then $P\,\cap\, \P^{n+2}$ is a point $x$ of $X$. Let $\pi :\P^{2n+1}\dasharrow \P^{n+2}$ be the projection $(x_i,y_j)\mapsto (x_i)$.  Since   the differentials of $Q_i$ and $q_i$ coincide at $x$,  the derivative $\pi_*$ maps $T_x(P)\subset T_x(\mathscr{X})$ into $T_x(X)$. Since $P$ is general, $\pi_*T_x(P)$ is a hyperplane in $T_x(X)$  -- this will follow from the proof of Proposition \ref{diagcomm} 1) below, where we construct explicitely pairs $(P,\lambda )$ with this property.
	
	 Therefore we have   a rational map
	\[\psi: \mathscr{F}\dasharrow \P T^*X\,\quad (P,\lambda )\mapsto (x=P\cap \P^{n+2}\, ,\  \pi_*T_x(P))\,.\]
	The symmetric group $\mathfrak{S}_{n-1}$ acts on $\mathbb{P}^{2n+1}$ by permuting the $y_j$, and the group $(\Z/2\Z)^{n-1}$ by changing their signs; this gives an action of the semi-direct product $G:=(\Z/2\Z)^{n-1}\rtimes \mathfrak{S}_{n-1}$. We make $G$ act on $\mathbb{A}^{n-1}$ through its quotient $\mathfrak{S}_{n-1}$, by permutation of the $\lambda _i$. This induces an action of $G$ on $\mathscr{X}$ and therefore on $\mathscr{F}$, compatible  via $p$ with the action on the base.
	The map  $\psi$  is   invariant under this   action, hence factors   through the  quotient $\mathscr{F}/G$.  By   passing to the  quotient we get a map
	$p^\sharp: \mathscr{F}/G\rightarrow  \mathbb{A}^{n-1} /\mathfrak{S}_{n-1}$.

	\begin{prop}\label{diagcomm}
		$1)\ \psi$ induces a birational map $\psi^{\sharp}:\mathscr{F}/G\dasharrow \P T^*X$.
		
		$2) $ There is a commutative diagram
		\[\xymatrix{\mathscr{F}/G \ar@{-->}[r]^{\psi^{\sharp}}\ar[d]_{p^{\sharp}}&\quad  \P T^*X\ar@{-->}@<2ex>[d]^{\varphi }\\ \mathbb{A}^{n-1}/\mathfrak{S}_{n-1}\ar[r]_{\sigma }^{\sim}&\mathbb{A}^{n-1}\subset  \P^{n-1}
		}\]where $p^{\sharp}$ is deduced from $p$, and $\sigma $ is the isomorphism given by symmetric functions.
	\end{prop}
	
	\pr 1) Let $(x,H)\in \P T^*X$; we want to describe  the pairs $(P,\lambda )$ such that $P\cap \P^{n+2}=\{x\} $ and $\pi_*T_x(P)=H$. The latter condition says that, via the decomposition
	\[T_x(\P^{2n+1})=T_x(\P^{n+2})\oplus \Ker \pi_*\,,\]
	$T_x(P)$ identifies with the graph of a linear map
	\[\alpha : H\rightarrow \Ker \pi _*\, .\]
	Using the basis $(\frac{\partial }{\partial y_1},\ldots ,\frac{\partial }{\partial y_{n-1}}  )$ of $\Ker \pi _*$, we have $\alpha =(\alpha _1,\ldots ,\alpha _{n-1})$, where the $\alpha _i$ are linear forms on $H$. The condition $P\subset \mathscr{X}_{\lambda }$ implies that the hessians $h_{Q_1}(x)$ and $h_{Q_2}(x)$ vanish on $T_x(P)$, which gives
	\begin{equation}\label{diag}
		h_{q_1}(x)_{\mid H}=-\sum_i \alpha_i^2\ ,\quad h_{q_2}(x)_{\mid H}=-\sum_i\lambda_i\alpha_i^2\,.\end{equation}
	This is a  simultaneous diagonalization of the quadratic forms $h_{q_1}(x)_{\mid H}$ and   $h_{q_2}(x)_{\mid H}$; when they are in general  position,  this determines the $\lambda_i$  up to permutation and the  $\alpha_i$ up to sign and permutation, which proves 1).
	
	\smallskip	
	2) Let $(P,\lambda )\in \mathscr{F}$, and let $(x,H):=\psi(P,\lambda )$. According to Proposition \ref{geom}, $\varphi (x,H)$ is given by the $(n-1)$-uple of  quadrics $ q\in \Pi$
	such that 
	the form $h_{q}(x)_{\mid H}$ is singular. Using $(\alpha _1,\ldots ,\alpha_{n-1})$ as coordinates on $H$, we see from (\ref{diag}) that this $(n-1)$-uple is given by $(\lambda _1,\ldots ,\lambda _{n-1})$, which proves 2).\qed

	\subsection{Proof of Proposition \ref{birat}.}

	Let $\lambda $ be a general element of $\A^{n-1}$. Let us denote by $\Gamma $ the subgroup $(\Z/2\Z)^{n-1}$ of $G$. From Proposition \ref{diagcomm} and the cartesian diagram
	\[\xymatrix{\mathscr{F}/\Gamma \ar[r]\ar[d]^p& \mathscr{F}/G\ar[d]^{p^{\sharp}}\\
		\A^{n-1}\ar[r] & \A^{n-1}/\mathfrak{S}_{n-1}
	}\]we see that the fiber $\varphi ^{-1}(\lambda )$ is birational to the quotient $\mathscr{F}_{\lambda}/\Gamma $. By (\ref{reid}) $\mathscr{F}_{\lambda }$ is isomorphic to $JC_{\mu ,\lambda }$, and 
	one can choose the isomorphism so that  $\Gamma $ acts on 
	$JC_{\mu ,\lambda }$ as $\{\pm 1_{J}\} \times \Gamma ^+$, where $\Gamma ^+$ is a group  of translations by $2$-torsion elements. This proves the Proposition.\qed

	\bigskip	
	\section{Fibers of $\Phi$}\label{fibersPhi}
	\subsection{Results}\label{res}
	We keep the settings of the previous section. Recall that our parameter $\lambda $ lives in \break $\A^{n-1}\subset \sy^{n-1}\Pi\cong \P^{n-1} $. For $\lambda $ in $\A^{n-1}$, we denote by $\tilde{\lambda } $ a lift of $\lambda $ in $\C^n$ for the 
	quotient map $\C^n\mo\{0\}\rightarrow \P^{n-1} $.
	\begin{thm}\label{thm}
		Assume that $X$ is general. For $\lambda \in \A^{n-1}$ general, the fiber $\Phi ^{-1}(\tilde{\lambda}  )$ is isomorphic to $A\mo Z$, where\,:
		
		$\bullet\ A$ is the abelian variety quotient of $JC_{\mu ,\lambda }$ by a $2$-torsion subgroup, isomorphic to $(\Z/2\Z)^{n-2}\,;$
		
		$\bullet\ Z$ is a closed subvariety of codimension $\geq 2$ in $A$. 
	\end{thm}
	
	\begin{cor}\label{lag}
		For \emph{every} smooth complete intersection of two quadrics $X\subset \P^{n+2}$, the fibration $\Phi :X\rightarrow \C^n$ is Lagrangian.
	\end{cor}
	\pr Assume first that $X$ is general. The symplectic form on $T^*X$ is $d\eta $, where $\eta $ is the Liouville form. 
	By the Theorem and the Hartogs principle, the pull back of $\eta $ to a general fiber of $\Phi $ is the restriction of a 1-form on an abelian variety, hence is closed. This implies the result.
	
	Let $p:\mathscr{X}\rightarrow B$ be a complete family of smooth intersection of two quadrics in $\P^{n+2}$. The constructions of \S \ref{Phi} can be globalized over $B$: we have a rank 2 vector bundle $\mathscr{V}$ over $B$ whose fiber at a point $b\in B$ is the space of quadratic forms vanishing on $\mathscr{X}_b$. We get a homomorphism $\sy^{n-1}\mathscr{V}\rightarrow p_*T_{\mathscr{X}/B}$, which gives rise to a morphism $\boldsymbol{\Phi}:T^*(\mathscr{X}/B)\rightarrow \sy^{n-1}\mathscr{V}^*$ over $B$ which induces over each point $b\in B$ our map $\Phi $. There is a natural Liouville form $\boldsymbol{\eta }$ on $T^*(\mathscr{X}/B)$; since $d\boldsymbol{\eta }$ vanishes on a general fiber of $\boldsymbol{\Phi}$, it vanishes on all fibers.\qed
	
	\begin{cor}
		Assume that $X$ is general. The multiplication map $\sy^{\pp} H^0(X,\sy^2T_X)\rightarrow H^0(X,\sy^{\pp}T_X)$ is an isomorphism.
	\end{cor}
	(We will give in \S\, \ref{tensors} a proof valid with no generality assumption.)
	\smallskip	
	
	\pr The Theorem implies that every function on a general fiber of $\Phi $ is constant, hence the pull back $\Phi ^*:H^0(\C^n,\O_{\C^n})\rightarrow H^0(T^*X,\O_{T^*X})$ is an isomorphism. The right hand space is canonically isomorphic to $H^0(X,\sy^{\pp}T_X)$, hence we get an algebra isomorphism  $\C[t_1,\ldots ,t_n]\iso H^0(X,\sy^{\pp}T_X)$. By construction the $t_i$ are mapped to elements of $H^0(X,\sy^{2}T_X)$, so the Corollary follows.\qed
	
	\begin{rem}\label{acis}
		Let $V_1,\ldots ,V_n$ be the Hamiltonian vector fields on $T^*X$ associated to the components of $\Phi $. For  $\lambda$ general in $\C^n$, let us identify $\Phi ^{-1}(\lambda )$ to $A\mo Z$ as in the Theorem. Then by Hartogs' principle the $V_i$ \emph{linearize} on $A$ --- that is, they extend to a basis  of $H^0(A,T_A)$. This allows in principle to write explicit solutions of the Hamilton equations for $\Phi _i$ in terms of theta function.
	\end{rem}

	\subsection{Proof of the Theorem: lemmas}
	We fix a general point $\lambda \in\A^{n-1}$. We denote by $\mathscr{F}\oo$ the open subset of $\mathscr{F}$ where the rational map $\psi$ is well-defined, and by $\mathscr{F}\oo_{\lambda }$ its intersection with the fiber $\mathscr{F}_{\lambda }$. Since $\lambda $ is general, the complement of $\mathscr{F}\oo_{\lambda }$ in $\mathscr{F}_{\lambda }$ has codimension $\geq 2$. The rational map $\psi$ induces a morphism $\psi\oo: \mathscr{F}\oo\rightarrow \P T^*X$; we denote by $\psi\oo_{\lambda }$ its restriction to $\mathscr{F}_{\lambda }\oo$. Let $Z\subset \P T^*X$ be the indeterminacy locus of $\varphi $ (\S \,\ref{Phi}), and let $\mathscr{F}_{\lambda }^{\operatorname{bad} }:= (\psi\oo_{\lambda })^{-1}(Z)\subset \mathscr{F}\oo_{\lambda }$. 
	\begin{prop}\label{codim}
		$\mathscr{F}_{\lambda }^{\operatorname{bad} }$ has codimension $\geq 2$ in $\mathscr{F}_{\lambda }$.
	\end{prop}
	We postpone the proof of  the Proposition to the next section, and first show how it implies Theorem \ref{thm}.
	
	Let $0_X\subset T^*X$ be the zero section, and let $q:T^*X\mo 0_X\rightarrow \P T^*X$ be the quotient map. Let $\varphi \oo: \P T^*X\mo Z\rightarrow \P^{n-1}$ be the morphism induced by $\varphi $. We have
	$\ q(\Phi^{-1}(\tilde{\lambda}   ))= (\varphi \oo)^{-1}(\lambda)$, and the restriction \[q_{\lambda }: \Phi^{-1}(\tilde{\lambda}   )\rightarrow  (\varphi \oo)^{-1}(\lambda)\] is an \'etale double cover, with Galois involution $\iota $ induced by $(-1_{T^*X})$.
	
	We put $\mathscr{F}_{\lambda }\ooo:= \mathscr{F}_{\lambda }\oo\mo \mathscr{F}_{\lambda }^{\operatorname{bad} }$, and consider the restriction \[\psi_{\lambda} \oo:\mathscr{F}_{\lambda }\ooo\rightarrow (\varphi \oo)^{-1}(\lambda )\quad\mbox{of } \ \psi\oo\,.\]
	
	\begin{lem}\label{lemlag}
		The fiber $\Phi ^{-1}(\tilde{\lambda})$ is Lagrangian, and has trivial tangent bundle.
	\end{lem}
	\pr The \'{e}tale double cover $q_{\lambda }$ induces by fibered product an \'{e}tale double cover
	\[\pi :\widetilde{\mathscr{F}}_{\lambda }\ooo\rightarrow \mathscr{F}_{\lambda }\ooo\]such that $\psi\oo_{\lambda }$ lifts to a morphism $\tilde{\psi}\oo_{\lambda }: \widetilde{\mathscr{F}}_{\lambda }\ooo\rightarrow \Phi^{-1}(\tilde{\lambda}   )$.  
	
	By Proposition \ref{codim}, the complement of $\mathscr{F}_{\lambda }\ooo$ in $\mathscr{F}_{\lambda }$ has codimension $\geq 2$, so $\pi $ extends to an \'etale double cover $\widetilde{\mathscr{F}}_{\lambda }\rightarrow \mathscr{F}_{\lambda }$, where $\widetilde{\mathscr{F}}_{\lambda }$
	is   an abelian variety or the disjoint union of two abelian varieties. The morphism $\tilde{\psi}\oo_{\lambda }: \widetilde{\mathscr{F}}_{\lambda }\ooo\rightarrow \Phi^{-1}(\tilde{\lambda}  )$ is generically of maximal rank. Again by Proposition  \ref{codim}, the holomorphic 1-forms on $\widetilde{\mathscr{F}}_{\lambda }\ooo$ are closed, hence by pull back the same holds for the holomorphic 1-forms on $\Phi^{-1}(\tilde{\lambda} )$. As in the proof of Corollary \ref{lag}, this implies that $\Phi^{-1}(\tilde{\lambda } )$ is  Lagrangian. The second assertion is a basic property of Lagrangian fibers.  \qed

	\begin{lem}\label{lift}
		The morphism $\psi\oo_{\lambda }$ lifts to a morphism $\tilde{\psi}\oo_{\lambda }:\mathscr{F}\ooo_{\lambda }\rightarrow \Phi ^{-1}(\tilde{\lambda}  ) $.
	\end{lem}
	
	\pr  It suffices to show that the double covering $\pi :\widetilde{\mathscr{F}}_{\lambda }\ooo\rightarrow \mathscr{F}_{\lambda }\ooo$ splits.
	
	Assume the contrary, so that $\widetilde{\mathscr{F}}_{\lambda }$ is an abelian variety.  By Lemma \ref{lemlag} $H^0(\Phi^{-1}(\tilde{\lambda}) ,\Omega ^1)$ has dimension $n$. It follows that the pull back 
	$(\tilde{\psi}\oo_{\lambda })^*: H^0(\Phi^{-1}(\tilde{\lambda}) ,\Omega ^1) \rightarrow H^0(\widetilde{\mathscr{F}}_{\lambda }\ooo,\Omega ^1)$ is bijective. Since the Galois involution of the double covering $\pi $ acts trivially on 
	holomorphic 1-forms, the same holds for
	the  Galois involution $\iota $ of the double covering $q_{\lambda }: \Phi^{-1}(\tilde{\lambda }  )\rightarrow  (\varphi \oo)^{-1}(\lambda)$. 
	
	Now we observe that the 1-forms on $\Phi^{-1}(\tilde{\lambda}) $ are ``pure", that is, extend to any smooth projective compactification of $\Phi^{-1}(\tilde{\lambda} )  $: this follows from the fact that this holds after pull back to $\widetilde{\mathscr{F}}_{\lambda }\ooo$. But the quotient $\Phi^{-1}(\tilde{\lambda } )/\iota $ is isomorphic to a Zariski open subset of $\varphi ^{-1}(\lambda )$, which by Proposition \ref{birat} has  no nonzero  holomorphic 1-forms, so that any Zariski open set has no nonzero closed pure holomorphic 1-forms.
	This contradiction proves the Lemma.\qed
	
	\subsection{Proof of  Theorem \ref{thm}}
	
	Lemma \ref{lift} gives a factorization   \[\psi\oo_{\lambda }: \mathscr{F}\ooo_{\lambda } \xrightarrow{\ \tilde{\psi}\oo_{\lambda } \ } \Phi ^{-1}(\tilde{\lambda}   )\xrightarrow{\ q_{\lambda }\ }  (\varphi \oo)^{-1}(\lambda)\, .\]   
	By Proposition \ref{birat}, $\psi\oo_{\lambda }$ induces a birational morphism
	\[\psi\oo_{\lambda ,\Gamma }: \mathscr{F}\ooo_{\lambda }/\Gamma \longrightarrow (\varphi \oo)^{-1}(\lambda )\,;\]
	it follows that for some subgroup $\Gamma '\subset \Gamma $ of index $2$, the morphism $\tilde{\psi_{\lambda}\oo}:\mathscr{F}_\lambda\ooo\rightarrow \Phi^{-1}(\tilde{\lambda} )$ factors through a birational morphism
	\[\tilde{\psi}\oo_{\lambda ,H'}: \mathscr{F}\ooo_{\lambda }/\Gamma '\longrightarrow \Phi ^{-1}(\tilde{\lambda} )\,.\]
	By Lemma \ref{lemlag},  the cotangent bundle of $\Phi ^{-1}(\tilde{\lambda } )$ is trivial. Therefore the cotangent bundle of $\mathscr{F}_{\lambda }\ooo/\Gamma '$ is generically generated by its global sections. This implies that $\Gamma '$ acts trivially on holomorphic 1-forms, hence is the subgroup $\Gamma ^{+}$ of $\Gamma $ generated by translations, isomorphic to $(\Z/2\Z)^{n-2}$; thus $\mathscr{F}_{\lambda }/\Gamma '$ is an abelian variety $A$. 
	
	To simplify notation, we put $A\oo:=\mathscr{F}\ooo_{\lambda }/\Gamma '$ and 
	$u:=\tilde{\psi}\oo_{\lambda ,H'}$. The rational map $u^{-1}: \Phi ^{-1}(\tilde{\lambda} )\dasharrow A$ is everywhere defined (see e.g. \cite[Theorem 4.9.4]{B-L}), so we have two morphisms
	\[A\oo \xrightarrow{\ u\ } \Phi ^{-1}(\tilde{\lambda} ) \xrightarrow{\ u^{-1}\ } A\]whose composition is the inclusion $A\oo\hookrightarrow A$.
	Since the tangent bundles of $A$ and $\Phi ^{-1}(\tilde{\lambda} )$ are trivial, the determinant of $Tu:  T_{A\oo}\rightarrow u^*T _{\Phi ^{-1}(\tilde{\lambda} )}$ is a function on $A\oo$, hence constant by Proposition \ref{codim}. Therefore $u$ is \'etale and birational, hence an open embedding. This implies that every function on $\Phi ^{-1}(\widetilde{\lambda})$ is constant (because its restriction to $A\oo$ is  constant). Then the previous argument shows that $u^{-1}$ is also an open embedding, so that   $\Phi ^{-1}(\tilde{\lambda} )$ is isomorphic to an open subset of $A$ containing $A\oo$. This proves the Theorem.

	\bigskip	
	\section{Proof of Proposition \ref{codim}}
	We keep the notations of (\ref{aux}).
	Recall that we have coordinates $(x_0,\ldots ,x_{n+2};y_1,\ldots ,y_{n-1})$ on $\P^{2n+1}$, and  subspaces $\P^{n+2}$ and $\P^{n-2}$ in $\P^{2n+1}$
	defined by $y=0$ and $x=0$.
	
	Let $q_1(x)=q_2(x)=0$ be the equations defining $X$ in $\P^{n+2}$, and
	let $R$ be the vector space of quadratic forms in $y=(y_1,\ldots ,y_{n-1})$. We define an extended family $\mathscr{X}^{e}\subset \P^{2n+1}\times R^2$ by  
	\[\mathscr{X}^{e}=\{\bigl((x,y);(r_1,r_2)\bigr)\in \P^{2n+1}\times R^2\ |\ 
	q_1(x)+r_1(y)=q_2(x)+r_2(y)=0\}\, .\]
	The fiber  $\mathscr{X}^{e}_{r}$ at a point $r=(r_1,r_2)$ of $R^2$ is the intersection in $\P^{2n+1}$ of the two quadrics $q_1(x)+r_1(y)=\allowbreak q_2(x)+r_2(y)=0$.
	Let $\G$ be the Grassmannian of $(n-1)$-planes in $\P^{2n+1}$;
	we define as before 
	\[\mathscr{F}^{e}:=\{(P,r)\in \G\times R^2\,|\,P \subset \mathscr{X}^{e}_{r}\}\]
	and the extended rational map $\psi^{e}:\mathscr{F}^{e}\dasharrow \P T^*X$, which maps a general $P\subset \mathscr{X}^{e}_r$ to the pair $(x,H)$ with $\{x\}=P\cap\P^{n+2} $, $H=\pi _*T_x(P)$.
	
	We observe that a general pair $r=(r_1,r_2)$ of $R^2$ is simultaneously diagonalizable, so the restriction of $\psi^{e}$ to $\mathscr{F}^{e}_r$ coincides, for an appropriate choice of the coordinates $(y_i)$, with the map $\psi_{\lambda }$ that we want to study. Thus Proposition \ref{codim} will follow from the following Proposition:
	\begin{prop}
		Assume that $X$ is general.
		
		$1)$ Let $\Gamma \subset \mathscr{F}^{e}$ be the locus of points $(P,r)$ such that either $\dim P\cap\P^{n+2}>0$, or $P\cap\P^{n-2}\neq \varnothing$. Then $\Gamma $ has codimension $\geq 2$ in $\mathscr{F}^{e}$.
		
		$2)$  There  exists no  divisor  in $\mathscr{F}^{e}\mo \Gamma$ which dominates $R^2$ and  is mapped to the base-locus $Z\subset \P T^*X$ by $\psi_e$.
	\end{prop}
	\pr 1) Let $Q$ be the vector space of quadratic forms on $\P^{2n+1}$ of the form $q(x)+r(y)$ for some quadratic forms $q$ and $r$. For each pair of integers $(k,l)$ with $k\geq 0$, $l\geq -1$, let $\G_{k,l}$ be the locally closed subvariety of $(n-1)$-planes $P\in\G$ such that
	\[\dim(P\cap \P^{n+2})=k\ ,\quad \dim (P\cap \P^{n-2})=l\,.\]
	(We put by convention $l=-1$ if $P\cap\P^{n-2}=\varnothing$.) Let
	\[ \mathscr{F}^{Q}:= \{(P, (Q_1,Q_2))\in \G\times Q^2\ |\ Q_{1|P}=Q_{2|P}=0\}\,, \]
	\[\mathscr{F}^{Q}_{k,l}:= \mathscr{F}^{Q}\cap (\G_{k,l}\times Q^2)\,.\]
	The general fiber of the  projection $\mathscr{F}^{Q}\rightarrow Q^2$  is an abelian variety, 
	and we recover $\mathscr{F}^{e}$ by restricting $\mathscr{F}^{Q}$ to pairs of quadratic forms of the form $(q_{1}(x)+r_{1}(y), q_{2}(x)+r_{2}(y))$. It thus suffices to prove the result for the larger family $\mathscr{F}^{Q}$, that is, to show that
	$\mathscr{F}_{k,l}^{Q}$ has codimension $\geq 2$ in $\mathscr{F}^{Q}$.

	This is done by a dimension count. For $P\in \G$, let $\varphi _P$ be the restriction map $Q\rightarrow H^0(P,\O_{P}(2))$. The fiber of the projection $\mathscr{F}^{Q}\rightarrow \G$ is the vector space $(\Ker \varphi _P)^{\oplus 2}$. For $P$ general, $\varphi _P$ is surjective: this is the case for instance if $P$ is contained in the  $(n+2)$-plane in $\P^{2n+1}$ defined by $y_i=x_i\ (i=1,\ldots ,n-1)$. However $\varphi _P$ is not surjective for $P\in \G_{k,l}$, because  the forms $r(y)_{|P}$  are singular along 
	$P\cap \mathbb{P}^{n+2}$ and the forms $q(x)_{|P}$ are singular along 
	$P\cap \mathbb{P}^{n-2}$: this implies that the subspaces $P\cap \mathbb{P}^{n+2}$ and $P\cap \mathbb{P}^{n-2}$ are apolar for all forms in $\im \varphi _P$. Therefore the corank of $\varphi _P$ is $\geq (k+1)(l+1)$, and there is equality when $P$ is contained in the subspace defined by $x_0=\ldots =x_{n+1-k}=y_1=\ldots =y_{n-2-l}=0$, hence for $P$ general in $\G_{k,l}$. Thus our assertion follows from:
	\begin{align*}
		\codim(\mathscr{F}^{Q}_{k,l},\mathscr{F}^{Q})&=\codim(\G_{k,l},\G)-2(k+1)(l+1) \\
		&= k(k+1)+(l+1)(l+4)-2(k+1)(l+1)\\
		&=(k-l)(k-l-1)+2(l+1)\\
		&\geq 2\quad \mbox{if }\ k\geq 1\ \mbox{ or }\ l\geq 0\,. 
	\end{align*}  
	
	\smallskip	
	2) The base locus $Z\subset \P T^*X$ has codimension $\geq 2$ (Corollary \ref{codimZ}). Note that $\psi^{e}$ is well-defined in $\mathscr{F}^{e}\mo\Gamma $. If $\mathscr{D}$ is a divisor in $\mathscr{F}^{e}\mo\Gamma $ with $\psi^{e}(\mathscr{D})\subset Z$, the map $\psi^{e}$ has not maximal rank along $\mathscr{D}$. This contradicts the following Lemma:
	
	\begin{lem}
		$\psi^{e}$ has maximal rank on $\mathscr{F}^{e}\mo\Gamma $.
	\end{lem}
	\pr Let $(x,H)$ be a point of $T^*X$; we view $H$ as a hyperplane in the projective tangent space to $x$ at $X$. The fiber of $\psi^{e}: \mathscr{F}^{e}\mo\Gamma \rightarrow \P T^*X$ at $(x,H)$ is the locus
	\makeatletter\tagsleft@false\makeatother
	\begin{align}
		(\psi^{e})^{-1}(x,H)=\{(P,r_1,r_2)\in \G\times R^2\ |\ P\cap \P^{n+2}&=\{x\}\, ,\  P\cap \P^{n-2}=\varnothing\, ,\ \pi (P)=H,    \label{1}
		\\ (q_i(x)+r_i(y))_{|P}&=0\quad (i=1,2)\} \,. \label{2}
	\end{align} 
	The equations (\ref{1}) define a smooth, locally closed subvariety $\G_{x,H}$ of $\G$. Let $P\in \G_{x,H}$, and let  $\chi _P: R\rightarrow H^0(P,\O_P(2))$ be the restriction map. We will show below that the image of $\chi  _P$ is the space of quadratic forms on $P$ which are singular at $x$. Since the forms $q_{i|P}$ are singular at $x$, this implies that the solutions of (\ref{2}) form an affine space over $(\Ker \chi_P)^{\oplus 2}$. Therefore $(\psi^{e})^{-1}(x,H)$ admits an affine fibration over $\G_{x,H}$, hence is smooth.
	
	Clearly the quadrics in $\im \chi  _P$ are singular at $x$. 
	To prove the opposite inclusion, choose the coordinates $(x_i)$ so that $x=(1,0,\ldots ,0)$. Since $P\cap \P^{n+2}=\{x\} $, there exist linear forms $\ell_1,\ldots ,\ell_{n+2}$ in the $y_j$ so that $P$ is defined by $x_i=\ell_{i}(y)$ for $i=1,\ldots ,n+2$. Then a quadratic form on $\P^{2n+1}$ singular at $x$ can be written as a form in $x_1,\ldots,x_{n+2} ;y_1,\ldots ,y_{n-1}$, hence its restriction to $P$ is in $\im \chi _P$. This proves the Lemma, hence also  the Proposition.\qed

	\bigskip
	\section{Symmetric tensors: second approach}\label{tensors}

\subsection{The cotangent bundle of a smooth quadric}\label{quadric}

We consider a smooth quadric $Q\subset \P^{n+1}$, defined by an equation $q=0$. Its cotangent bundle $\P T^*Q$  parameterizes pairs $(x,P)$ with $x\in Q$ and $P$ a $(n-1)$-plane tangent to $Q$ at $x$. Thus we get a morphism $\gamma $ from $\P T^*Q$ to the grassmannian $\G$ of $(n-1)$-planes in $\P^{n+1}$, which is the morphism defined by the linear system $\abs{\O_{\P T^*Q}(1)}$. It is birational onto its image, but contracts the subvariety $\mathscr{C}\subset \P T^*Q$ consisting of pairs $(x,P)$ such that $P$ is tangent to $Q$ along a line $\ell\subset Q$, and $x\in\ell$: then   $\gamma ^{-1}(P)$ consists of the pairs $(x,P)$ with $x\in\ell$.

Let $h_q\in H^0(Q,\sy^2 \Omega^1 _Q(2))$ be the hessian form of $q$ (\S \ref{Phi}). 
Choosing coordinates $(x_i)$ such that $q(x)=\sum x_i^2$, we have $h_q=\sum (dx_i)^2$ (note that this is, up to a scalar, the unique element of $H^0(Q,\sy^2 \Omega^1 _Q(2))$ invariant under $\Aut(Q)$). Then  $h_q(x)$ is non-degenerate at each point $x$ of $Q$, so 
 $h_q$ induces an isomorphism $\Omega^1 _Q(1)\iso T_Q(-1)$, hence also $\sy^2 \Omega^1 _Q(2)\iso \sy^2 T_Q(-2)$. The image in $H^0(Q,\sy^2 T_Q(-2))$ of $h_q$ by this isomorphism is $h'_q=\sum \partial _j^2$. We will view $h'_q$ as an element of $H^0(\P T^*Q,\O_{\P T^*Q}(2)\otimes p^*\O_Q(-2))$, where $p: \P T^*Q\rightarrow Q$ is the projection.
\begin{prop}\label{divhq}
The divisor of $h'_q$ is $\mathscr{C}$.  The projection $p_{|\mathscr{C}}:\mathscr{C}\rightarrow Q$ is a smooth quadric fibration, and 
 $\mathscr{C}$ is a prime divisor for $n\geq 3$.
\end{prop}

\pr Let $x\in Q$; the hyperplane tangent to $x$ at $Q$ cuts down a cone over the smooth quadric $Q_x\subset \P(T_x(Q))$ defined by $h_q(x)=0$ (\S\,\ref{Phi}). The isomorphism $T_x(Q)\iso T^*_x(Q)$ given by $h_q(x)$ carries $Q_x$ into the dual quadric $Q_x^*$ in   $\P(T^*_x(Q))$. On the other hand, a point  $y\in p^{-1}(x)$ corresponds to a hyperplane  $H_y\subset \P(T_x(Q))$,   and $y$
 belongs to $\mathscr{C}$ if and only if $H_y$ is tangent to $Q_x$, that is $y\in Q_x^*$. This proves the equality $\mathscr{C}=\operatorname{div} (h'_q) $. Thus the fiber of  
$p_{|\mathscr{C}}:\mathscr{C}\rightarrow Q$ at $x$ is $Q_x$, which is smooth, and connected if $n\geq 3$.\qed

\begin{rem}
The variety $\mathscr{C}$ is an example of a total dual VMRT \cite{H-L-S}, for the proof of the Theorem we will combine this tool with the birational transformation of $\P T^* X$ defined by a double cover, cf.\ \cite{A-H}.
\end{rem}
\medskip	
We will have to consider the following situation. Let $Q'$ be another quadric in $\P^{n+1}$, such that the intersection $B:=Q\cap Q'$ is a
smooth hypersurface in $Q$. The surjection $T_{Q}\rightarrow N_{B/Q}$ gives a section of $\P T^*Q$ over $B$, hence an embedding $s:B\hookrightarrow \P T^*Q$.
\begin{lem}\label{BC}
The image $s(B)$ is not contained in $\mathscr{C}$.
\end{lem}
\pr Let $x\in B$. The point $s(x)$ in $\P(T_x^*(Q))$ corresponds to the hyperplane 
image of $T_x(B)$ in $T_{x}(Q)$; we must show that this hyperplane is not tangent to the quadric $Q_x:=h_q(x)$. In terms of projective space, this means that the projective tangent space to $Q'$ at $x$ is not tangent to the cone $Q\cap \P T_x(Q)$ at a smooth point $y$ of $Q$.

Suppose this is the case, with $y=(y_0,\ldots ,y_{n+1})$. We can assume that $Q'$ is defined by $\sum \alpha _i x_i^2=0$, with $\alpha _i\in\C$ distinct. Then the (projective) tangent space to $Q'$ at $x$, given by $\sum (\alpha _ix_i)\xi _i=0$, must coincide with the tangent space  to $Q$ at $y$, given by $\sum y_i\xi _i=0$. This implies $y=(\alpha _0x_0,\ldots ,\alpha _{n+1}x_{n+1})$. Thus the point $x$ must satisfy
\[\sum x_i^{2}=\sum \alpha _ix_i^2=\sum \alpha _i^2 x_i^2=0\,.\]
If these relations hold for all $x$ in $B$, the quadric $\sum \alpha _i^2 x_i^2=0$ must belong to the pencil spanned by $Q$ and $Q'$. This  means that there exist scalars $\lambda ,\mu ,\nu$ such that 
\[\lambda \alpha _i^2+\mu \alpha _i+\nu=0\quad \mbox{for all }i\,,\]
which is impossible since the $\alpha _i$ are distinct. Therefore there exists $x\in B$ such that $s(x)\notin \mathscr{C}$.\qed

\subsection{Explicit description of symmetric tensors}
	\label{ss.explicit-description}

We keep the notation of the previous sections: $X\subset \P=\P^{n+2}$ is defined by $q_1=q_2=0$, with
\[
	q_1 = \sum_{i=0}^{n+2} x_i^2, \qquad
	q_2 = \sum_{i=0}^{n+2} \mu_i x_i^2 
	\qquad
	\mbox{with } \mu_i \in \C \mbox{ distinct}.
	\]
We put $\partial _i:=\dfrac{\partial }{\partial x_i} $. We have an exact sequence
\[0\rightarrow T_X\rightarrow T_{\P|X}\xrightarrow{\ (dq_1,dq_2)\ }\O_X(2)^2\rightarrow 0\,,\]where $dq_i$ maps the restriction of a vector field $V$ on $\P$ to $V\cdot q_i$. This gives an exact sequence of symmetric tensors
\begin{equation}\label{sy2}
0\rightarrow \sy^2T_X\rightarrow \sy^2 T_{\P|X}\xrightarrow{\ (dq_1,dq_2)\ }T_{\P|X}(2)^2\,,
\end{equation} 
where $dq_i(V_1V_2)=(V_1\cdot q_i)V_2+(V_2\cdot q_i)V_1$ for $V_1,V_2$ in $H^0(X, T_{\P|X})$.
\begin{prop}\label{basis}
 The quadratic vector fields $\displaystyle s_i:=\sum_{j\neq i} \dfrac{(x_i\partial _j-x_j\partial _i)^2}{\mu _j-\mu _i} $ in $H^0(X, \sy^2 T_{\P|X})$ belong to the image of $H^0(X, \sy^2 T_{X})$.
\end{prop}
\pr  According to the exact sequence (\ref{sy2}) we have to prove $dq_1(s_i)=dq_2(s_i)=0$.

We have $(x_i\partial _j-x_j\partial _i)\cdot q_1=0$, hence $dq_1(s_i)=0$, and   $(x_i\partial _j-x_j\partial _i)^2\cdot q_2=2(\mu _j-\mu _i)x_ix_j(x_i\partial _j-x_j\partial _i)$, hence, using $\sum x_j\partial _j=0$ and $q_{1|X}=0$:
\[ dq_2(s_i)=2x_i^2\sum_{j\neq i} x_j\partial _j-(2x_i\partial _i)\sum_{j\neq i}x_j^2=0\,, \quad\mbox{which proves the Proposition.}\eqno{\qed}\]
\smallskip	
Fron now on we will consider the $s_i$ as elements of $H^0(X, \sy^2 T_{X})$.

\medskip	
 \subsection{The double cover}\label{pi}
 
Let $p: \P^{n+2}\dasharrow \P^{n+1}$ be the projection $(x_0,\ldots ,x_{n+2})\mapsto (x_1,\ldots ,x_{n+2})$. The image $p(X)$ is the smooth quadric $Q$ in $\P^{n+1}$ defined by 
\[\sum_{i=1}^{n+2}(\mu _i-\mu _0)x_{i}^{2}=0\,.\]
The restriction $\pi : X\rightarrow Q$ of $p$ is a double covering,  branched along the subvariety $B\subset Q$ defined by
\[\sum_{i=1}^{n+2} x_i^2=\sum_{i=1}^{n+2}\mu _i x_i^2=0\,.\]
It is a smooth complete intersection of 2 quadrics in $\P^{n+1}$. The ramification locus $R\subset X$ of $\pi $ (isomorphic to $B$)  is the hyperplane section $x_0=0$ of $X$.

The tangent map of $\pi :X\rightarrow Q$ gives a morphism
\[ \tau : T_X\rightarrow \pi ^*T_Q\]
which is an isomorphism outside of $R$. 
Consider the normal exact sequence
\[0\rightarrow T_R \rightarrow T_{X|R}\rightarrow N_{R/X}\rightarrow 0\,.\]
The involution $\iota:(x_0,\ldots ,x_{n+2})\mapsto (-x_0,x_1,\ldots ,x_{n+2}) $  acts on $T_{X|R}$; this splits the exact sequence, giving a decomposition
\[T_{X|R}=T_R\oplus N_{R/X}\]into eigenspaces for the eigenvalues $+1$ and $-1$. Let $\rho : T_{X|R}\rightarrow T_R$ be the projection on the first summand. We deduce from $\rho $ a sequence of homomorphisms
\[h^k: H^0(X, \sy^kT_X)\longrightarrow H^0(X, \sy^kT_{X|R})\xrightarrow{\ \sy^k\rho \ } H^0(R, \sy^kT_R)\,.\]
Since $\iota _*\partial _0=-\partial _0$ and $\iota _*\partial _j=\partial _j$ for $j>0$, we have 
\begin{equation}\label{restriction}
h^2(s_0)=0\quad\mbox{and}\quad h^2(s_i)=\sum_{j>0 \atop j\neq i} \dfrac{(x_i\partial _j-x_j\partial _i)^2}{\mu _j-\mu _i}\quad\mbox{for }i>0\,;
\end{equation}
in other words, $h^2$ maps $s_1,\ldots ,s_{n+2}$   to the elements $\hat{s}_1,\ldots ,\hat{s}_{n+2}$ of $H^0(R, \sy^2T_R)$ constructed in Proposition \ref{basis} applied to $R$.
\medskip	

Let $\pi ^*\P T^*Q$ be the pull back under $\pi $ of the projective bundle $\P T^*Q\rightarrow Q$.
The homomorphism $\tau :T_X\rightarrow \pi ^*T_Q$ gives rise to a birational map $g: \pi ^*\P T^*Q\dasharrow \P T^*X$.
Following the geometric description of the tangent map as an elementary transformations of vector bundles in the sense of Maruyama \cite{M1},\cite[Corollary 1.1.1]{M2}, one has a commutative diagram
	\begin{equation} \label{diagram-transform}
		\xymatrix{
			& \Gamma
			\ar[dl]_{\mu}  \ar[dr]^{\nu}
			&  
			\\
			\pi ^*\P T^*Q\ar@{-->}[rr]^{g}
			\ar[rd]_{p}
			& 
			& \mathbb{P}T^* X\ar[ld]^{q}
			\\
			& X &   
		}    
	\end{equation}
	where $p$ and $q$ are the canonical projections, $\nu: \Gamma\rightarrow \P T^* X$ is the blow-up along the subspace
	$\P T^* R \subset \P T^* X$ defined by the projection $\rho $, $\mu :\Gamma\rightarrow \pi ^*\P T^*Q$
	is the blow-up of the image $B'$ of the embedding $B\hookrightarrow \pi ^*\P T^*Q$ deduced from the surjective homomorphism $\pi ^*T_Q\rightarrow \pi ^*N_{B/X}$.
	
	Let $E_{\mu}$ be the exceptional divisor of $\mu$. By \cite[Theorem 1.1]{M2}, there is an isomorphism
	\begin{equation}
		\label{transform-tautological}
		\mu^*\O_{\pi ^*\P T^*Q}(1) \otimes \O_{\Gamma}(-E_{\mu}) \cong \nu^*\O_{\P T^* X}(1)\, ,
			\end{equation}
	as well as the equality
	\begin{equation}
		\label{transform-exceptional}
		\nu_* E_{\mu}
		=
		q^*R\, .
	\end{equation}

\smallskip	
\subsection{The divisor of $s_0$}
We now consider the divisor $\mathscr{C}\subset \P T^*Q$ defined in (\ref{quadric}), and the cartesian diagram
\[\xymatrix{\pi ^*\P T^*Q\ar[r]^{\pi '}\ar[d] & \P T^*Q\ar[d]\\X \ar[r]^{\pi }&Q \,.}\]
Put 
$\mathscr{C}':= \pi '^{-1}(\mathscr{C})$.  The projection $\mathscr{C}'\rightarrow X$ is again a smooth quadric fibration, so $\mathscr{C}'$ is smooth, and connected for $n\geq 3$.

Recall that we have defined the element $\displaystyle s_0:=\sum_{j=1}^{n+2} \dfrac{(x_0\partial _j-x_j\partial _0)^2}{\mu _j-\mu _0} \in H^0(X,\sy^2T_X)$ (\ref{ss.explicit-description}). We will view $s_0$ as an element of $H^0(\P T^*X,\O(2))$.
\begin{prop}\label{p.construct-Di}
Assume $n\geq 3$. We have $g_*\mathscr{C}'=\divi(s_0)$. 
\end{prop}
\pr We first show that $g_*\mathscr{C}'\in\abs{\O_{\P T^*X}(2)}$. By Proposition \ref{divhq} we have $\mathscr{C}'\in \abs{\O_{\pi ^*\P T^*Q}(2)\otimes p^*\O_X(-2)}$. Using (\ref{transform-tautological}), (\ref{transform-exceptional}) and the projection formula, we get the linear equivalences
\[\nu_*\mu ^*\mathscr{C}' \sim  2\nu_*\mu ^*(c_1(\O_{\pi ^*\P T^*Q}(1)-p^*R))\sim 2(c_1(\O_{\P T^*X}(1)) +q^*R)-2q^*R=c_1(\O_{\P T^*X}(2))\,.\]
Thus it is enough to prove that $\nu_*\mu ^*\mathscr{C}' $ is irreducible.  Since $\mathscr{C}'$ is irreducible and $\mu$ is the blow-up along $B'\subset \pi ^*\P T^*Q$, it suffices to show that  $B'$ is not contained in $\mathscr{C}'$. If this is the case, we have $\pi '(B')\subset \pi' (\mathscr{C}')=\mathscr{C}$. But $\pi '(B')=s(B)$, where $s: B\hookrightarrow \P T^*Q$ is the embedding defined by the surjective homomorphism $T_Q\rightarrow N_{B/Q}$. Then the result follows from Lemma \ref{BC}.

Since $g_*\mathscr{C}'$ and $\divi(s_0)$ are linearly equivalent effective  divisors and $g_*\mathscr{C}'$ is irreducible, it suffices to show that their restrictions to $\P T_x^*X$ coincide for a general point $x\in X$.
		
		Fix a point $x=[x_0,\dots,x_{n+2}]\in X\mo R$, so that $x_0\not =0$.
	Then the tangent map $T\pi (x):T_x(X)\rightarrow T_{\pi (x)}(Q)$ is an isomorphism;  in the diagram (\ref{diagram-transform}), the maps $\mu ,\nu$ and $g$ restricted over the fibers at $x$ are all isomorphisms. Let us show that $\mathscr{C}'$ and $T\pi (\divi(s_0))$ define the same quadric in $\P (T_{\pi (x)}(Q))$.
	
	Now $\mathscr{C}'\cap \P (T^*_{x}(X))=\mathscr{C}\cap \P (T^*_{\pi (x)}(Q))$ is  the  quadric defined by the element $h'_q$ of (\ref{quadric}). In the coordinates $(z_i)$ defined by $z_i= (\mu _i-\mu _0)^{1/2}x_i$, the equation of $Q$ is $\displaystyle\sum_{j=1}^{n+2}z_j^2=0$, so
	\[h'_q=\sum_{j=1}^{n+2}\bigl(\dfrac{\partial }{\partial z_j}\bigr)^2= \sum_{j=1}^{n+2}\dfrac{\partial _j^2}{\mu _j-\mu _0}\ \cdot \]
On the other hand, since $\pi (x_0,\ldots ,x_{n+2})=(x_1,\ldots ,x_{n+2})$, we have $T\pi (\partial _0)=0$ and $T\pi (\partial _j)=\partial _j$ for $j>0$, hence
\[T\pi (s_0)=x_0^2\, \sum_{j=1}^{n+2}\dfrac{\partial _j^2}{\mu _j-\mu _0}\ \cdot  \]	
Since $x_0\neq 0$, this proves the Proposition.\qed

\subsection{Proof of part a) of the Theorem}

	Suppose now that $n \geq 3$. Consider the double cover
	$
	\pi  : X \rightarrow Q$
	and the ramification divisor  $R \subset X$. The restriction maps $h^k$ defined in  \eqref{pi} yield a homomorphism of graded $\C$-algebras 
	\[
	h: S(  X) := H^0(X, \sy^{\pp} T_X) 
	\longrightarrow H^0(R, \sy^{\pp} T_{R}) =: S( R).
	\]
	\begin{prop}
The kernel $\mathscr{I}$ of $h$
	is the ideal generated by $s_0$. 
\end{prop}
\pr Since  $\mathscr{I}$ is a homogeneous ideal, it suffices  to prove that every homogeneous element $s\in \mathscr{I}$ can be written as $s=s' s_0$ for some element $s'\in S(X)$. 

Fix an element $s \in \mathscr{I}$ of degree $k$. It corresponds to an effective Cartier divisor $G$ in  the linear system
	$|\O_{\P T^* X}(k)|$. Recall the commutative diagram \eqref{diagram-transform}  
\begin{equation*} 		\xymatrix{
			& \Gamma
			\ar[dl]_{\mu}  \ar[dr]^{\nu}
			&  
			\\
			\pi ^*\P T^*Q\ar@{-->}[rr]^{g}
			\ar[rd]_{p}
			& 
			& \mathbb{P}T^* X	\ar[ld]^{q}
			\\
			& X &   
		}    
	\end{equation*}
Put $ \hat{G} :=\mu_* \nu^* G\subset \pi ^*\P T^*Q$.
 By \eqref{transform-tautological}, $\hat{G}$ belongs to the linear system
	$\abs{\O_{\pi ^*\P T^*Q}(k)}$.

	Here comes the key observation: since $s \in \mathscr{I}$, the divisor
	$ \hat{G} \subset \pi ^*\P T^*Q$
	contains $p^* R$. Indeed, 
	since $(\pi ^*T_Q)_{|R}$ is invariant under $\iota $, the homomorphism $\tau^{} _{|R}$ factors as
\[\tau^{} _{|R}:T_{X|R}\xrightarrow{\ \rho \ }T_R \longrightarrow  (\pi ^*T_Q)_{|R}\,.\]
Therefore	we have a commutative diagram
	\[\xymatrix{H^0(X,\sy^kT_X)\ar[r]^{h^k}\ar[d]_{\sy^k\tau }& H^0(R,\sy^kT_R)\ar[d]\\
	H^0(X,\sy^k\pi ^*T_Q)\ar[r]& H^0(R,\sy^k (\pi ^*T_Q)_{|R})
	}\]so that $\sy^k\tau (s)$ vanishes on $R$. But $\hat{G}$ is the divisor of $\sy^k\tau (s)$, viewed as a section of $\O_{\pi ^*\P T^*Q}(k)$, hence
 $\hat{G}$ contains $p^*R$.
 
 \medskip	
 
 Now we want to show that the divisor $\mathscr{C}'\subset \pi ^*\P T^*Q$ is a component of $\hat{G}-p^*R$. Recall (\ref{quadric}) that $\mathscr{C}$ is the union of the lines $\ell$ which are contracted by the morphism $\gamma : \P T^*Q\rightarrow \G$, so that $c_1(\O_{\P T^*Q}(1))\cdot \ell=0$. Thus the curves $\ell':=\pi '^*\ell$ cover $\mathscr{C}'$, and satisfy $c_1(\O_{\pi ^*\P T^*Q}(1))\cdot \ell'=0$. On the other hand the divisor $R\subset X$ is a hyperplane section, so $p^*R\cdot \ell'=R\cdot p_*\ell'>0$. Therefore
 \[(\hat{G}-p^*R)\cdot \ell'<0\, ,\]so $\mathscr{C}'$ is a component of $\hat{G}$. Thus $g_*\mathscr{C}'$ is a component of $G$. Since $g_*\mathscr{C}'=\divi(s_0)$ by Proposition \ref{p.construct-Di}, this proves  the Proposition.\qed
 
 \medskip	
 The following Proposition implies part a) of our main Theorem: 
  \begin{prop}
Assume $n\geq 2$. For any choice of indices $0\leq i_1<\ldots <i_n\leq n+2$,
the homomorphism $\C[t_1,\ldots ,t_{n}]\rightarrow S(X)$ which maps $t_j$ to $s_{i_{j}}$,  with $\deg(t_i)=2$, is an isomorphism of graded $\C$-algebras.
\end{prop}

\pr We argue by induction on $n$. The statement for $n=2$ follows from \cite[Theorem 5.1]{DO-L}, except the fact that any two of the $s_i$ generate $H^0(X,\sy^2T_X)$. Up to permuting of the coordinates, it suffices to prove that $s_0$ and $s_1$ are linearly independent. But $h^2:H^0(X,\sy^2T_X)\rightarrow H^0(R,\sy^2T_R)$ maps $s_0$ to zero and $s_i$, for $i>0$, to the corresponding elements
$\hat{s}_i$ of  $H^0(R,\sy^2T_R)$; this implies our assertion.

Assume $n\geq 3$. 
By the induction hypothesis,  the homomorphism $\C[t_1,\ldots ,t_{n-1}]\rightarrow S(R)$ which maps $t_i$ to $\hat{s_i}$ is an isomorphism of graded $\C$-algebras (with $\deg(t_i)=2$). 
It follows that $h$ is surjective, and that $(s_0,\ldots ,s_{n-1})$   form a basis of $H^0(X,\sy^2 T_X)$ and generate the $\C$-algebra $S(X)$. Thus we have a surjective homomorphism  $u:\C[t_0,\ldots ,t_{n-1}]\rightarrow S(X)$, with $u(t_i)=s_i$.
	 
	  In particular, the Krull dimension of $S(X)$ is at most $n$. On the other hand,  the ring $S(X)$ is a domain and $s_0$ is neither zero nor a unit. Thus, by Krull's Hauptidealsatz, the Krull dimension of $S(X)$ is equal to $n$, hence $u$ is an isomorphism. By permutation of the coordinates we get the same result for any choice of $n$ elements in $\{s_0,\ldots ,s_{n+2}\} $, hence the Proposition.\qed


\bigskip	
	\begin{thebibliography}{X-X-X}
		
		\bibitem[A-H]{A-H} F. Anella, A. H\"oring\,: \textsl{The cotangent bundle of a K3 surface of degree two}. Preprint \texttt{arXiv:2207.09294}.

		\bibitem[B-H-K]{B-H-K}  I. Biswas, Y.I. Holla, C. Kumar\,: \textsl{On moduli spaces of parabolic vector bundles of rank 2 over} $\C \P^1$. 
		Michigan Math. J. \textbf{59}  (2010), no. 2, 467-475. 
		
		\bibitem[B-L]{B-L} C. Birkenhake, H. Lange\,: \textsl{Complex abelian varieties}. 
		Grundlehren der mathematischen Wissenschaften \textbf{302}. Springer-Verlag, Berlin, 1992. 
		
		\bibitem[B-N-R]{BNR} A. Beauville, M.S. Narasimhan, S. Ramanan\,: \textsl{Spectral curves and the generalised theta divisor}. J. Reine Angew. Math. \textbf{398}  (1989), 169-179.
		
		
		\bibitem[C]{C} C. Casagrande\,: \textsl{Rank $2$ quasiparabolic vector bundles on $\P^1$ and the variety of linear subspaces contained in two odd-dimensional quadrics}. Math. Z. \textbf{280}  (2015), no. 3-4, 981-988.
		
		
		\bibitem[DO-L]{DO-L} B. De  Oliveira, C. Langdon\,: \textsl{Twisted symmetric differentials and the quadric algebra of subvarieties of $\P^N$ of low codimension}. Eur. J. Math. \textbf{5}  (2019), no. 2, 454-475. 
		
		\bibitem[D-R]{D-R}  U. Desale, S. Ramanan\,: \textsl{Classification of vector bundles of rank $2$ on hyperelliptic curves}. Invent. Math. \textbf{38}  (1976/77), no. 2, 161-185. 
		
		
		
		\bibitem[G-H]{G-H} P. Griffiths,J. Harris\,: \textsl{Algebraic geometry and local differential geometry}. Ann. Sci. \'Ecole Norm. Sup. (4) \textbf{12}  (1979), no. 3, 355-452.
		
		\bibitem[H]{H} N. Hitchin\,: \textsl{
			Stable bundles and integrable systems}. 
		Duke Math. J. \textbf{54}  (1987), no. 1, 91-114.
		
		\bibitem[H-L-S]{H-L-S} A. H\"oring, J. Liu, F. Shao\,: \textsl{Examples of Fano manifolds with non-pseudoeffective tangent bundle}. J. Lond. Math. Soc. (2) \textbf{106}  (2022), no. 1, 27-59.
		
		\bibitem[K-L]{K-L} H. Kim, Y. Lee\,: \textsl{Lagrangian fibration structure on the cotangent bundle of a del Pezzo surface of degree $4$}. Preprint \texttt{arXiv:2210.01317}.

\bibitem[M1]{M1}
M. Maruyama\,: \textsl{On a family of algebraic vector bundles}, Dissertation Kyoto University,
  \url{https://doi.org/10.14989/doctor.r2072}, 1972.
  
		\bibitem[M2]{M2} M. Maruyama\,: \textsl{On a family of algebraic vector bundles}. Number theory, algebraic geometry and commutative algebra, in honor of {Y}asuo {A}kizuki, pp.~95-146.  Kinokuniya, Tokyo, 1973.
		  
  \bibitem[N]{N} P. Newstead\,:  \textsl{Stable bundles of rank $2$ and odd degree over a curve of genus} $ 2$. 
		Topology \textbf{7}  (1968), 205-215. 
		
	
		
		\bibitem[R]{R} M. Reid\,: \textsl{The complete intersection of two or more quadrics}. PhD thesis, Cambridge University, 1972.
		
				
		\bibitem[V]{V} P. Vanhaecke\,: \textsl{Integrable systems in the realm of algebraic geometry}. Lecture Notes in Mathematics \textbf{1638}. Springer-Verlag, Berlin, 1996.
		
		
	
	\end{thebibliography}
\end{document}